\documentclass[amssymb]{conm-p-l}
\begin{document}
\newcommand{\bref}[1]{$\mbox{(\ref{#1})}$}
\newcommand{\p}{{\Bbb {R^+}} \rightarrow {\Bbb {C}}}
\newcommand{\q}{{\mathcal Q}}
\newcommand{\z}{{\Bbb {Z}} ^d}
\newcommand{\s}{\sum\limits_{j =1}^{d}}
\newtheorem{thm}{Theorem}
\newtheorem{lem}{Lemma}
\newtheorem{conj}{Conjecture}
\newtheorem{rem}{Remark}
\newtheorem{cor}{Corollary}

\title{The Uncertainty Principle for relatively dense sets and lacunary spectra}

\author{Oleg Kovrizhkin}

\address{MIT, Math Dept, 2-273, 77 Mass Ave, Cambridge, MA 02139, USA}

\email{oleg@math.mit.edu}


\date{}


\begin{abstract} We obtain a new version of the Uncertainty Principle for functions with Fourier transforms supported on a lacunary set
of intervals. This is a generalization of Zygmund's theorem on lacunary
trigonometric series to the real line in the spirit of the Logvinenko-Sereda theorem for relatively dense sets.
\end{abstract}

\maketitle

\section{Introduction}\label{intro}
We prove a new version of the Uncertainty Principle which is a general statement saying that a function and its Fourier transform can not be simultaneously 
concentrated on small sets. Many examples of the Uncertainty Principle can be found in the text book by Havin and J\"oricke \cite{HJ} and a paper by Folland 
and 
Sitaram \cite{FS}. We restrict 
ourselves to the following type of the Uncertainty Principle:
\begin{eqnarray}
\int_E |f|^2 \ge C\|f\|_2^2\label{upineq}
\end{eqnarray}
for all $f$ with supp$\hat f \subset \Sigma$ where $E^c$ and $\Sigma$ are small subsets of the real line and $C = C(E, \Sigma) > 0$ does not depend on $f$.
Our result is a combination of two versions of this principle. One is Zygmund's theorem on lacunary trigonometric series (\cite{Z}, pp. 202-208) and the other 
is the 
Logvinenko-Sereda
theorem for relatively dense sets (\cite{HJ}, p.113), \cite{LS}. As other examples of the inequality \bref{upineq}, we mention the Amrein-Berthier theorem 
\cite{AB}, (\cite{HJ}, p. 97, p. 455), \cite{N} (which is a quantitative 
version of a result due to Benedicks \cite{B}) where $E^c$ and $\Sigma$ are sets of finite measure  and Wolff's theorem \cite{W} where $E^c$ and $\Sigma$ are 
so called $\epsilon$-thin sets.

A sequence of integers $\Lambda = \{n_i\}_{i = -\infty}^{\infty} \subset {\Bbb {Z}}$ is called lacunary with parameter $R$ if 
$$\sup\limits_{r \neq 0} Card\{(i,j): n_i - n_j = r\} = R < +\infty.$$ For example, a sequence of integers $\Lambda$ satisfying  the Hadamard condition $\frac 
{n_{i+1}} {n_i} \ge q > 1$ is lacunary (\cite{Z}, p. 203).

~

{\bf Zygmund's Theorem:} {\it Given any lacunary sequence $\Lambda = \{n_i\}_{-\infty}^{\infty}$ and a measurable set $E \subset [0,1]$ of positive measure 
then for any $g(x) \in L^2[0,1]$ with spec $g \subset \Lambda$: $g(x) = \sum c_{n_j} e^{i2\pi n_j x}$ we have
$$\int_E |g|^2 \ge C(E, R) \|g\|_2^2,$$
where $C(E, \Lambda) > 0$ depends only on $E$ and $\Lambda$.}

Moreover, Nazarov showed \cite{N} that actually $C(E, \Lambda)$ can be replaced with $C(|E|, R)$.

~

A measurable set $E \subset {\Bbb {R}}$ is called relatively dense if there are $\gamma >0$ and $a > 0$ such that
$$|E \cap I| \ge \gamma |I|$$
for any interval $I$ of length $a$.

~

{\bf Logvinenko-Sereda Theorem:} {\it Given any relatively dense set $E \subset {\Bbb {R}}$ and any function $f \in L^2({\Bbb{R}})$ such that supp $\hat f 
\subset [-b, 
b]$ then
$$\int_E |f|^2 \ge C(\gamma, a , b)\|f\|_2^2$$
where $C(\gamma, a , b) > 0 $ depends only on $\gamma$, $a$ and $b$.
}

In his recent paper \cite {OK1} the author obtained an optimal estimate of $C(\gamma, a , b) = \left ( \frac {\gamma}{C} \right )^{C(ab + 1)}$ and generalized 
the theorem to functions whose Fourier transform 
is supported in a union of finitely many intervals with an estimate on $C$ depending only on the number of intervals but not how they are placed. The results 
also hold in higher dimensions.

Our main results here are the following two theorems:

\begin{thm} If   $\Lambda = \{n_i\}_{i = -\infty}^{\infty} \subset {\Bbb {Z}}$ is a lacunary sequence with parameter $R$ then there exist $\epsilon (R) < 1$ 
and 
an 
absolute constant $C > 0$
such that for any $f \in L^2$ with supp $\hat f \subset \bigcup\limits_{i} [n_i - b/4, n_i + b/4]$ and for any relatively dense set $E \subset {\Bbb {R}}$ 
which satisfies $|E^c \cap I| \le \epsilon (R) |I|$ for any interval $I$ of length $1/b$ provided $b \le 1$  or $|E^c \cap I| \le \frac {\epsilon (R)} b |I|$ 
for any interval $I$ of length $1/b$ provided $b > 1$ we have
$$\int_E |f|^2 \ge C\|f\|_2^2$$
where $C > 0$ is an absolute constant.
\end{thm}

\begin{thm}
Let $E$ be a periodic set with period 1 such that $|E \cap [0,1]| = \gamma > 0$. If $\Lambda = \{n_i\}_{i = -\infty}^{\infty} \subset {\Bbb {Z}}$ is a 
lacunary 
sequence with parameter $R$ and $f \in L^2$ with supp $\hat f \subset \bigcup\limits_{i} [n_i - 1/2, n_i + 1/2]$ then
$$\int_E |f|^2 \ge C(\gamma, R)\|f\|_2^2$$
where $C(\gamma, R) > 0$ depends only on $\gamma$ and $R$.
\end{thm}

The constant $C$ below is not fixed and might change appropriately from one equality or inequality to another one.

\section{Proof of Theorem 1}\label{proof1}
First we will fix some notations. Let $g(x)$ denote a 1-periodic $L^2[0,1]$ function with lacunary spectrum: spec $g \subset \Lambda$, i.e., $g(x) = \sum 
c_{n_i}e^{i 2\pi n_i x}$ where $c_n$ stands for the $n$-th Fourier coefficient of $g$. Let $\phi$ be a fixed $C_{0}^{\infty}$ function with supp $\phi \subset 
[-1/2, 1/2]$ and such that $\phi(x) \equiv 1$ when $x \in 
[-1/4, 1/4]$. It is clear that $|\check {\phi} (x)| \le \frac C {1 + x^2}$. Denote $\phi _n (x) = e^{i2\pi n x} \phi (x)$. Therefore $|\check {\phi_n} (x)| =  
|\check {\phi} (x + n)| \le \frac C {1 + (x + n)^2}$.

\begin{lem}
$$\|g\|_4 \le (1 + R)^{\frac 1 4} \|g\|_2$$ 
\end{lem}

{\bf Proof:} This is true since
\begin{eqnarray*}
\int_0^1 |g|^4 &=& \int_0^1 (|g|^2)^2 = \sum\limits_n |\widehat{|g|^2}(n)|^2 \\
&=& \sum\limits_n |\sum\limits_{k-l = n}c_k \bar c_l|^2 = |\sum\limits_{k}|c_k|^2|^2 + \sum\limits_{n \neq 0} |\sum\limits_{k-l = n}c_k \bar c_l|^2.
\end{eqnarray*}
Using H\"older's inequality we can estimate the second term by
$\sum\limits_{n \neq 0} \sum\limits_{k-l = n}R|c_k|^2 |c_l|^2 \le R\sum\limits_{k} \sum\limits_{l}|c_k|^2 |c_l|^2 = R(\sum\limits_{k}|c_k|^2)^2$. Thus,

$$\int_0^1 |g|^4 \le (\sum\limits_{k}|c_k|^2)^2 + R(\sum\limits_{k}|c_k|^2)^2 = (1 + R)\|g\|_2^4.$$
\hfill$\square$

\begin{lem}
Let $E \subset {\Bbb {R}}$ be a relatively dense set which satisfies $|E^c \cap I| \le \epsilon |I|$ for any interval $I$ of length $1/b$ provided $b \le 1$  
or $|E^c \cap I| \le \frac {\epsilon } b |I|$ for any interval $I$ of length $1/b$ provided $b > 1$ then
\begin{eqnarray}
b\int_{E^c} |g(x)|^2 \cdot |\check {\phi_k} (bx)|^2 dx \le C(R) \sqrt{\epsilon}\|g\|_2^2\label{l2}
\end{eqnarray}
where $C(R)$ depends only on $R$. 
\end{lem}

{\bf Proof:} We have
\begin{eqnarray}
b\int_{E^c} |g(x)|^2 \cdot |\check {\phi_k} (bx)|^2 dx &=& \sum\limits_n b\int\limits_{E^c \cap [n/b, (n+1)/b]}|g(x)|^2 \cdot |\check {\phi_k} (bx)|^2 
dx\nonumber\\
&\le& b \sum\limits_n \sqrt{\int\limits_{[n/b, (n+1)/b]}|g(x)|^4 dx} \cdot \sqrt{\int\limits_{E^c \cap [n/b, (n+1)/b]}|\check {\phi_k} (bx)|^4 dx}\nonumber\\
&\le& b \sum\limits_n \sqrt{(1 + 1/b)\|g\|_4^4} \cdot \sqrt{1/b \int\limits_{b\cdot E^c \cap [n, n+1]}|\check {\phi_k} (x)|^4 dx} \label{l21}.
\end{eqnarray}
We used H\"older's inequality to obtain the first inequality. In the second inequality we covered the interval $[n/b, (n+1)/b$ by no more than $1 + 1/b$ 
intervals of length 1.  Using {\bf Lemma 1} and the fact that $|\check {\phi_k} (x)|  \le \frac C {1 + (x + k)^2}$ we can estimate \bref{l21} by
\begin{eqnarray}
&\le& C(R) \sqrt{b + 1}\cdot\|g\|_2^2 \sum\limits_n \sqrt{\int\limits_{b\cdot E^c \cap [n, n+1]}\frac C {(1 + (x + k)^2)^2} dx}\nonumber\\
&\le& C(R) \sqrt{b + 1}\cdot\|g\|_2^2 \sum\limits_n \frac 1 {1 + (n + k)^2}\sqrt{|b\cdot E^c \cap [n, n+1]|}\label{l22}.
\end{eqnarray}
Note that $|b\cdot E^c \cap [n, n+1]| \le \epsilon$ if $b \le 1$  and $|b\cdot E^c \cap [n, n+1]| \le \epsilon / b$ if $b > 1$. In both cases we can bound 
\bref{l22} by
\begin{eqnarray*}&\le& C(R)\|g\|_2^2 \sum\limits_n \frac 1 {1 + (n + k)^2}\cdot \sqrt{\epsilon}\\
&=& C(R)\|g\|_2^2 \sum\limits_n \frac 1 {1 + n^2}\cdot \sqrt{\epsilon} = C(R)\sqrt{\epsilon}\|g\|_2^2. 
\end{eqnarray*}
\hfill$\square$

\begin{lem}
Let $E \subset {\Bbb {R}}$ be a relatively dense set which satisfies $|E^c \cap I| \le \epsilon |I|$ for any interval $I$ of length $1/b$ provided $b \le 1$  
or $|E^c \cap I| \le \frac {\epsilon } b |I|$ for any interval $I$ of length $1/b$ provided $b > 1$ then for $k \neq l$
\begin{eqnarray}
b\int_{E^c} |g(x)|^2 \cdot |\check {\phi_k} (bx) \cdot \check {\phi_l} (bx)| dx \le C(R) \sqrt{\epsilon}\|g\|_2^2 \cdot \frac 1 {|k - l|^2}\label{l3}
\end{eqnarray}
where $C(R)$ depends only on $R$.
\end{lem}

The proof is similar to the one of {\bf Lemma 2}. Just use the facts that $|\check {\phi_k} (x)| \le \frac C {1 + (x + k)^2}$ and 
\begin{eqnarray*}
\sum\limits_n \frac 1 {(1 + (n+k)^2)(1 + (n+l)^2)} &=& \sum\limits_n \frac 1 {(1 + n^2)(1 + (n+l - k)^2)}\\
&\le& C\int\frac 1 {(1 + x^2)(1 + (x+l - k)^2)}dx \le \frac C {|k - l|^2}.
\end{eqnarray*}

Now we are in a position to proceed with the proof of {\bf Theorem 1}. Since supp$\hat f \subset \bigcup\limits_{n_i} [n_i - b/4, n_i + b/4]$ we can choose 
(not necessarily uniquely)

$\hat f_{n_i} \in L^2$ such that supp$\hat f_{n_i} \subset [n_i - b/4, n_i + b/4]$, the supports of $\hat f_{n_i}$ are disjoint and $\hat f = 
\sum\limits_{n_i}\hat 
f_{n_i}$. Note that $\hat f_{n_i} = \hat f \cdot \chi_{[n_i - b/4, n_i + b/4]}$ if $|i|$ is large enough since the intervals $[n_i - b/4, n_i + b/4]$ and 
$[n_j - 
b/4, n_j + b/4]$ are disjoint for $i \neq j$ if $|i|$ is large enough, which follows from the fact that there are no more than $R\cdot b$ pairs $(i,j)$ such 
that $0 < |n_i - n_j| \le b/2$. Although $\hat f_{n_i}$ is supported on $[n_i - b/4, n_i + b/4]$, we will define its Fourier series on the larger interval 
$[n_i - 
b/2, 
n_i + b/2]$ converging in $L^2$ as follows:

$$\hat f_{n_i}(x) = \sum\limits_k c_{n_i}^{(k)} e^{\frac {i2\pi k x} b} \cdot \chi_{[n_i - b/2, n_i + b/2]}(x)$$ where $c_{n_i}^{(k)} = \frac 1 b 
\int\limits_{n_i - \frac b 2}^{n_i + \frac b 2} \hat f_{n_i}(x) e^{\frac{-i2\pi k x}{b}}$.

Therefore,

$$\|f_{n_i}\|_2^2 = \|\hat f_{n_i}\|_2^2 = b\sum\limits_k |c_{n_i}^{(k)}|^2.$$

Since the supports of $\hat f_{n_i}$ are disjoint, we have
\begin{eqnarray*}
\|f\|_2^2 &=& \|\hat f\|_2^2 = \sum\limits_{n_i} \|\hat f_{n_i}\|_2^2\\
&=& \sum\limits_{n_i} b\sum\limits_k |c_{n_i}^{(k)}|^2 = b \sum\limits_k \|g_k\|_2^2
\end{eqnarray*}

where $g_k(x) = \sum\limits_{n_i} c_{n_i}^{(k)}e^{i2\pi n_i x}$ are 1-periodic functions with lacunary spectra: spec $g_k \subset \Lambda$ and $\|g_k\|_2^2 = 
\int\limits_0^1 |g_k(x)|^2 dx$.

Using the facts that $\phi(\frac {x - n_i} b) \equiv 1$ when $x \in [n_i - b/4, n_i + b/4]$ and supp $\hat f_{n_i} \subset [n_i - b/4, n_i + b/4]$, we get
\begin{eqnarray*}\hat f_{n_i}(x) &=& \hat f_{n_i}(x) \cdot \phi(\frac {x - n_i} b)\\
&=& \sum\limits_k c_{n_i}^{(k)} e^{\frac {i2\pi k x} b} 
\cdot \chi_{[n_i - b/2, n_i + b/2]}\cdot \phi(\frac {x - n_i} b).
\end{eqnarray*}

Now use the fact that supp $\phi(\frac {\cdot - n_i} b) \subset [n_i - b/2, n_i + b/2]$ to obtain from the previous equality the following:
\begin{eqnarray*}
\hat f_{n_i}(x) 
&=& \sum\limits_k c_{n_i}^{(k)} e^{\frac {i2\pi k x} b} \cdot \phi(\frac {x - n_i} b)\\
&=& \sum\limits_k c_{n_i}^{(k)}\cdot \phi_k(\frac {x - n_i} b).
\end{eqnarray*}

Taking the inverse Fourier transform, we get
$$f_{n_i}(x) = b\sum\limits_k c_{n_i}^{(k)}e^{i2\pi n_i x}\check \phi_k(bx).$$

Therefore,
\begin{eqnarray*}
f(x) &=& \sum\limits_{n_i} f_{n_i}(x)\\
&=& \sum\limits_{n_i} (b\sum\limits_k c_{n_i}^{(k)}e^{i2\pi n_i x}\check{\phi_k}(bx))\\
&=& b \sum\limits_k \check{\phi_k}(bx)(\sum\limits_{n_i}c_{n_i}^{(k)}e^{i2\pi n_i x})\\
&=& b \sum\limits_k g_k (x) \check{\phi_k}(bx).
\end{eqnarray*}

Now we can estimate $\int_{E^c}|f|^2$:
\begin{eqnarray}
\int\limits_{E^c}|f|^2 &=& b^2 \int\limits_{E^c}|\sum\limits_k g_k (x) \check \phi_k(bx)|^2 dx\nonumber\\
&=& b^2 \int\limits_{E^c}(\sum\limits_k  |g_k (x) \check \phi_k(bx)|^2 + \sum\limits_{k \neq l}g_k (x) \check \phi_k(bx) \bar g_l (x) \bar{\check 
{\phi_l}}(bx))dx\label{t1}.
\end{eqnarray}
Using {\bf Lemma 2} we estimate the first term in \bref{t1}:
\begin{eqnarray}
b^2 \sum\limits_k \int_{E^c} |g_k (x) \check \phi_k(bx)|^2 dx &\le& C(R)\sqrt{\epsilon}b\sum\limits_k\|g_k\|_2^2\nonumber \\
&=& C(R)\sqrt{\epsilon}\|f\|_2^2\label{t15}.
\end{eqnarray}
Using {\bf Lemma 3} we estimate the second term in \bref{t1}:
\begin{eqnarray}
&&b^2 \sum\limits_{k \neq l} \int_{E^c}g_k (x) \check \phi_k(bx) \bar g_l (x) \bar{\check {\phi_l}}(bx)dx \le \nonumber \\ &&b^2 \frac 1 2 \sum\limits_{k \neq 
l}\int_{E^c}
( |g_k 
(x) |^2 + |g_l (x) |^2)\cdot |\check \phi_k(bx) {\check \phi_l}(bx)|dx \le \nonumber \\
&& b\cdot C(R)\sqrt{\epsilon}b\sum\limits_{k \neq l}\frac{(\|g_k\|_2^2 + \|g_l\|_2^2)}{|k - l|^2} = \nonumber \\
&& b\cdot C(R)\sqrt{\epsilon}\sum\limits_{n \neq 0}\sum\limits_{k} \frac{(\|g_k\|_2^2 + \|g_{k+n}\|_2^2)}{n^2} = \nonumber \\
&& C(R)\sqrt{\epsilon}b\sum\limits_{k} \|g_k\|_2^2 \cdot \sum\limits_{n
\neq 0} \frac 1 {n^2} \le \nonumber \\
&& C(R)\sqrt{\epsilon}\|f\|_2^2\label{t16}.
\end{eqnarray}
Adding the estimates \bref{t15} and \bref{t16} we get
$$\int_{E^c}|f|^2 \le C(R)\sqrt{\epsilon}\|f\|_2^2.$$
Now we choose $\epsilon$ such that $C(R)\sqrt{\epsilon} \le 1/2$. Hence,
$$\int_{E}|f|^2 \ge \frac 1 2 \|f\|_2^2.$$
\hfill$\square$

\section{Proof of Theorem 2}\label{proof2}
A similar case was studied for uniqueness
in \cite{LF}: if $f$ vanishes on $E$ then $f$ vanishes on the whole real line.\\
We will start with some results on periodizations. Define a family of
periodizations of a function $f \in L^1$:
\begin{eqnarray}g_t(x) = \sum\limits_{k = -
\infty}^{\infty}f(x + k) e ^{-i2\pi t(x + k)}\label{Period}\end{eqnarray}
where $t \in [-\frac 1 {2}, \frac 1 {2}]$. Then $g_t(x)$ is periodic with
period 1 , $\|g_t\|_{L^1([0,1])} \le \|f\|_1$ and its Fourier coefficients
are:
\begin{eqnarray*}\hat g_t(l) = \hat f ( l  +
t).\end{eqnarray*} Now we assume that $f \in L^1 \cap L^2$. The next argument shows an
important relation between the average of the $L^2$ norm of periodizations and the $L^2$
norm of $f$:
\begin{eqnarray*}\int\limits_{-\frac 1 {2}}^{\frac 1
{2}}\int\limits_{E \cap [0,1]}|g_t(x)|^2 dx dt &=&
\sum\limits_{k,l}\int\limits_{E \cap [0,1]}\int\limits_{-\frac 1 {2}}^{\frac 1
{2}}f(x + k) \bar f(x + l) e^{-i2\pi t(k-l)}dt dx\nonumber \\
&=&\sum\limits_{k,l}\delta_{k l}\int\limits_{E \cap [0,1]}f(x + k) \bar f(x +
l)dx\nonumber \\
&=&\sum\limits_{k = -
\infty}^{\infty}\int\limits_{(E - k) \cap [0,1]}|f(x + k)|^2 dx \nonumber \\
&=& \int\limits_{E}|f|^2 .\end{eqnarray*} We used that $E = E - k$ since $E$ is
1-periodic. In particular, it follows that
\begin{eqnarray*}\int\limits_{-\frac 1 {2}}^{\frac 1
{2}}\int\limits_{[0,1]}|g_t(x)|^2 dx dt &=& \int |f|^2
 .\end{eqnarray*}
In the next lemma we extend these results to functions from $L^2$. 
\begin{lem} If $f \in L^2$ then there exists a family $\{g_t(x)\}_{t \in [-\frac 1 2, \frac 1 2]}$ of periodic
functions: $x \in [0,1]$, a.e. $t \in [-\frac 12, \frac 12]$ with period $1$
such that $g_t(x) \in L^2([0,1] \times
[-\frac 12, \frac 12])$,
$$\int\limits_{\frac 12}^{\frac 12}\int\limits_{E \cap [0,1]}|g_t(x)|^2 dx dt =
\int\limits_{E}
|f|^2$$
and
$$\hat g_t(l) = \hat f (l + t)$$
for almost all $t$ and all $l \in \Bbb{Z}$.
\end{lem}

{\bf Proof:} Consider the cutoff:
$$f^{n}(x) = \chi_{[-n,n]}f(x).$$
Since $f^n \in L^1 \cap L^2$ and converge to $f$ in $L^2$ we can define
corresponding families of periodizations $g_t^n(x)$ which form a Cauchy sequence
in $L^2([0,1] \times
[-\frac 12, \frac 12])$:
$$\int\limits_{-\frac 12}^{\frac 12}\int\limits_{[0,1]}|g_t^n(x) - g_t^m(x)|^2
dx dt = \int |f^n - f^m|^2.$$ Let $g_t(x)$ be the $L^2([0,1] \times [-frac 12, frac 12])$ limit of $g_t^n(x)$. Thus, we get the
first statement of the lemma
\begin{eqnarray*}\int\limits_{-\frac 12}^{\frac 12}\int\limits_{E \cap
[0,1]}|g_t(x)|^2 dx dt &=& \lim\limits_{n \to \infty}
\int\limits_{-\frac 12}^{\frac 12}\int\limits_{E \cap [0,1]}|g_t^n(x)|^2 dx dt\\
&=& \lim\limits_{n \to \infty}\int\limits_{E} |f^n|^2 = \int\limits_{E}
|f|^2.\end{eqnarray*} To obtain the second statement of {\bf Lemma 4} we consider the following
sum:
\begin{eqnarray*}\sum\limits_{l = -
\infty}^{\infty}\int\limits_{-\frac 12}^{\frac 12}|\hat g_t(l) - \hat f(l +
t)|^2 dt
&=& \sum\limits_{l = -
\infty}^{\infty}\int\limits_{-\frac 12}^{\frac 12}|\hat g_t(l) - \hat g_t^n(l) +
\hat f^n
(l + t) - \hat f(l + t)|^2 dt\\
&\le& 2\sum\limits_{l = -
\infty}^{\infty}\int\limits_{-\frac 12}^{\frac 12}|\hat g_t(l) - \hat
g_t^n(l)|^2 + |\hat f^n
(l + t) - \hat f(l + t)|^2 dt\\
&=& 2\int\limits_{-\frac 12}^{\frac 12}\int\limits_{[0,1]}|g_t(x) - g_t^n(x)|^2 dxdt +
2\int |f^n - f|^2 \le \epsilon\end{eqnarray*} where $\epsilon$ can be arbitrarily small
if $n$ is large enough. \hfill$\square$

Let $g_t$ be a family of periodizations of $f$ as defined in
\bref{Period}. It follows from {\bf Lemma 4} that $\hat g_t(l) = \hat f(l + t)$ for almost all
$t \in [-\frac 12, \frac 12]$. Since supp$\hat f \subset \bigcup\limits_{n_i}[n_i-1/2, n_i + 1/2]$, $g_t$ has lacunary spectra for almost all $t$. Therefore, 
from Nazarov's result \cite{N} for Zygmund's theorem it follows that
\begin{eqnarray}\int\limits_{E \cap [0,1]}|g_t|^2 \ge C(\gamma,
R)\int\limits_{0}^{1}|g_t|^2\label{Almost}
\end{eqnarray}
for almost all $t \in 
[-\frac 12, \frac 12]$. Applying {\bf Lemma 4} and \bref{Almost} we obtain
\begin{eqnarray*}\int\limits_{E}|f|^2 \ge C(\gamma,
R)\int|f|^2.\end{eqnarray*}\hfill$\square$

~

\end{document}